\documentclass[12pt]{amsart}
\usepackage{amsmath, amssymb, latexsym, amsthm, mathrsfs, color}
\usepackage[top=2.5cm, bottom=2.5cm, left=2.5cm, right=2.5cm]{geometry}
\usepackage[all]{xy}

      \newenvironment{changemargin}[2]{\begin{list}{}{
         \setlength{\topsep}{0pt}\setlength{\leftmargin}{0pt}
         \setlength{\rightmargin}{0pt}
         \setlength{\listparindent}{\parindent}
         \setlength{\itemindent}{\parindent}
         \setlength{\parsep}{0pt plus 1pt}
         \addtolength{\leftmargin}{#1}\addtolength{\rightmargin}{#2}
         }\item }{\end{list}}

\newcommand{\nc}{\newcommand}

\nc{\thusfar}{\par\bigskip\centerline{\my{--- Edited thus far ---}}\par\bigskip}
\nc{\lei}{\le^\oo}
\nc{\card}[1]{\left|#1\right|}
\nc{\medcard}[1]{\biggl|\,#1\,\biggr|}
\nc{\smallcard}[1]{|\,#1\,|}
\nc{\bds}{bidirectional $\roth$-scale}
\nc{\PF}{\mathrm{PF}}
\nc{\bbT}{\mathbb{T}}
\nc{\bbN}{\mathbb{N}}%{\w}
\nc{\beq}{\begin{eqnarray*}}\nc{\eeq}{\end{eqnarray*}}
\nc{\mbq}{\mb{?}}
\nc{\mb}[1]{{\mbox{\textbf{#1}}}}
\nc{\nop}{$\times$}
\nc{\fbn}{\!\!\fbox{\!\nop\!}\!\!}
\nc{\yup}{\checkmark}
\nc{\forces}{\Vdash}
\nc{\name}[1]{\dot{#1}}
\nc{\tf}{\my{FINISHED THUS FAR}}
\nc{\FU}{Fr\'echet--Urysohn}
\nc{\gs}{$\gamma$~space}
\nc{\Ga}{\Gamma}\nc{\Om}{\Omega}
\nc{\smallbinom}[2]{\begin{psmallmatrix} #1\\ #2 \end{psmallmatrix}}
\nc{\bgamma}{\smallbinom{\Om}{\Ga}}
\newcommand{\two}{\{0,1\}}
\nc{\productive}[2]{\bigl(#1,\allowbreak #2\bigr)^\x}
\nc{\prdct}[1]{\bigl(#1\bigr)^\x}
\nc{\Sel}{\mathsf{S}}
\nc{\sset}[2]{\{\,#1 : #2\,\}}
\nc{\smb}[1]{{\!\!\mb{#1}\!\!}}
\nc{\medset}[2]{{\biggl\{\,#1 : #2\,\biggr\}}}
\nc{\smallmedset}[2]{{\bigl\{\,#1 : #2\,\bigr\}}}
\nc{\set}[2]{{\left\{\,#1 : #2\,\right\}}}
\nc{\sseq}[1]{\{\, #1 : n\in\bbN\}}
\nc{\seq}[2]{{\la\, #1 : #2\,\ra}}
\nc{\eseq}[1]{#1_1, \allowbreak #1_2, \allowbreak\dotsc} %explicit sequence
\nc{\cube}{(\Cantor)^\bbN}
\nc{\Match}{\op{Match}}
\nc{\concat}[1]{\hat{\phantom{a}}\langle #1\rangle}
\nc{\poset}{\mathbb{P}}
\nc{\fn}[1]{{\op{Fn}(#1\times\w,2)}}
\nc{\linadd}{\op{linadd}}
\nc{\nonprod}{\non^\x}
\nc{\alephes}{{\aleph_0}}
\nc{\my}[1]{\marginpar{\textcolor{red}{***}}\textcolor{red}{#1}}
\nc{\later}[1]{{\color{green} #1}}
\nc{\BTs}[1]{{\color{green} #1 (BT)}}
\nc{\Cp}{\op{C}_\mathrm{p}}
\nc{\Bp}{\op{B}_p}
\nc{\Pa}[8]{\bibitem{#1} {#2}, \emph{#3}, {#4} \textbf{#5} ({#6}), {#7}--{#8}.}
\nc{\tPa}[5]{\bibitem{#1} {#2}, \emph{#3}, {#4}, to appear.}
\nc{\sPa}[4]{\bibitem{#1} {#2}, \emph{#3}, {#4}, submitted.}
\nc{\Bc}[9]{\bibitem{#1} {#2}, \emph{#3}, in: \textbf{#4} (#5), #6 #7, #8--#9.}
\nc{\fD}{\mathfrak{D}}
\nc{\fX}{\mathfrak{X}}
\nc{\Onbd}{\Op_{\mathrm{nbd}}} %{\Op_{\mathsf{nbd}}}
\nc{\Omnb}{\Om_{\mathrm{nbd}}} %{\Om_{\mathsf{nbd}}}
\nc{\od}{\mathfrak{od}}
\nc{\Setting}[7]{\xymatrix@R=4pt@C=7pt{#1\ar@{-}[r]&#2\ar@{-}[r]&#3\\&#4\ar@{-}[u]\\
#5\ar@{-}[uu]\ar@{-}[r] & #6\ar@{-}[u]\ar@{-}[r] & #7\ar@{-}[uu]}}
\nc{\mx}[1]{\begin{matrix}#1\end{matrix}}
\nc{\plim}{p\txt{-}\lim}
\nc{\Bgp}{{\Z^\bbN}}
\nc{\Cgp}{{{\Z_2}^\bbN}}
\nc{\Cite}[1]{\textbf{[#1]}}
\nc{\Next}[1]{{#1^+}}
\nc{\cFin}{\mathrm{cF}}
\nc{\scsp}{\text{-scale space}}
\nc{\cfn}{\text{cofinal}\ }
\nc{\Con}{\text{Concentrated}}
\nc{\Lind}{\text{Lindel\"of}\,}
\nc{\con}{\text{-Concentrated}}
\nc{\lind}{\text{-Lindel\"of}\,}
\nc{\ctbl}{\text{countably}\,}
\nc{\Men}{\text{Menger}}
\nc{\men}{\text{-Menger}}
\nc{\Hur}{\text{Hurewicz}}
\nc{\intvl}[2]{{[#1(#2),\allowbreak #1(#2\!+\!1))}}
\nc{\Dfin}{\mathfrak{D}_\mathrm{fin}}
\nc{\grbl}{{\mathrm{\tiny gp}}}
\nc{\bbP}{\mathbb{P}}
\nc{\BOfat}{\B_{\Om_{\mathrm{fat}}}}%\B_{\mathrm{fat}}}
\nc{\Bgood}{\B_{\mathrm{good}}}
\nc{\compactN}{\cl{\mathbb{N}}}
\nc{\blocks}[2]{\op{cl}_{#2}(#1)}
\nc{\blocksplus}[2]{\op{cl}^+_{#2}(#1)}
\nc{\arx}[1]{\texttt{http://arxiv.org/math/#1}}
\nc{\bq}{\begin{quote}}
\nc{\eq}{\end{quote}}
\nc{\cl}[1]{\overline{#1}}
\nc{\CH}{the Continuum Hypothesis}
\nc{\MA}{Martin's Axiom}
\nc{\Bfat}{\B_\mathrm{fat}}
\nc{\inv}{^{-1}}
\nc{\Cantor}{{\two^\bbN}}%{{2^\w}}
\nc{\bP}{\mathbf{P}}
\nc{\bof}{\op{\fb}}
\nc{\dof}{\op{\fd}}
\nc{\bofF}{\bof(\cF)}
\nc{\sr}[3]{\underset{\text{#2}}{\mbox{#1}}}
\nc{\gp}{\binom{\Om}{\Ga}}
\nc{\gpsmall}{\mbox{$\gp$}}
\nc{\gig}{\gimel}%{\gimel\Ga}
\nc{\gns}{\sone(\Om,\gig)}
\nc{\nsr}[2]{#1}
\nc{\Srg}{{\mathbb{S}}}
\nc{\Srgs}{{\mathbb{S}^*}}
\nc{\NN}{{\bbN^{\bbN}}}
\nc{\ZN}{{\Z^{\bbN}}}
\nc{\NNup}{{\bbN^{\uparrow\bbN}}}
\nc{\Pof}{\op{P}}
\nc{\PN}{{\Pof(\bbN)}}
\nc{\roth}{{[\bbN]^{\mbox{\tiny $\infty$}}}} %{{[\w]^{\w}}}
\nc{\Fin}{[\bbN]^{\text{$<\!\!\infty$}}} %{{[\w]^{<\w}}}%{{[\bbN]^{<\aleph_0}}}
\nc{\ici}{[\bbN]^{ \infty, \infty}}%{{[\bbN]^{(\aleph_0,\aleph_0)}}}
\nc{\Inc}{{\compactN^{\uparrow\bbN}}}
\nc{\powInc}[1]{{\big(\Inc\big)^{#1}}}
\nc{\powFin}[1]{{\big(\Fin\big)^{#1}}}
\nc{\powPN}[1]{{\big(\PN\big)^{#1}}}
\nc{\NcompactN}{{\compactN^\bbN}}
\nc{\Uarrow}{\smash{\big\uparrow}}
\nc{\LE}{\preccurlyeq}
\nc{\GE}{\succcurlyeq}
\nc{\op}{\operatorname}
\nc{\im}{\op{im}}
\nc{\Span}{\op{span}}
\nc{\maxfin}{\op{maxfin}}
\nc{\ran}{\op{range}}
\nc{\iso}{\cong}
\nc{\Madd}{{\M}^\star}
\nc{\cI}{\mathcal{I}}
\nc{\cJ}{\mathcal{J}}
\nc{\scrA}{\mathscr{A}}
\nc{\scrB}{\mathscr{B}}
\nc{\scrC}{\mathscr{C}}
\nc{\scrD}{\mathscr{D}}
\nc{\scrF}{\mathscr{F}}
\nc{\scrK}{\mathscr{K}}
\nc{\A}{\forall}
\nc{\B}{\mathrm{B}}
\nc{\cB}{\mathcal{B}}
\nc{\bB}{\mathbf{B}}
\nc{\BS}{\mathbf{B}(\mathcal{S})}
\nc{\BF}{\mathbf{B}(\mathcal{F})}
\nc{\BU}{\mathbf{B}(\mathcal{U})}
\nc{\cSp}{\mathcal{S}^+}
\nc{\cFp}{\mathcal{F}^+}
\nc{\cUp}{\mathcal{U}^+}

\nc{\BG}{\B_\Ga}
\nc{\BL}{\B_\Lambda}
\nc{\BT}{\B_\Tau}
\nc{\BTstar}{\B_{\Tau^*}}
\nc{\BO}{\B_\Om}
\nc{\DO}{\cD_\Om}
\nc{\KO}{\cK_\Om}
\nc{\CG}{\rmC_\Ga}
\nc{\FG}{\rmF_\Ga}
\nc{\CL}{\rmC_\Lambda}
\nc{\CT}{\rmC_\Tau}
\nc{\CTstar}{C_{\Tau^*}}
\nc{\CO}{C_\Om}
\nc{\COgp}{C_{\Om^{\grbl}}}
\nc{\CLgp}{C_{\Lambda^{\grbl}}}
\nc{\BOgp}{\B_{\Om}^{\grbl}}
\nc{\BLgp}{\B_{\Lambda^{\grbl}}}
\nc{\sfC}{\mathsf{C}}
\nc{\sfD}{\mathsf{D}}
\nc{\bD}{\mathbf{D}}
\nc{\Tau}{\mathrm{T}}
\nc{\cA}{\mathcal{A}}
\nc{\cK}{\mathcal{K}}
\nc{\cD}{\mathcal{D}}
\nc{\cF}{\mathcal{F}}
\nc{\cS}{\mathcal{S}}
\nc{\cT}{\mathcal{T}}
\nc{\cG}{\mathcal{G}}
\nc{\cY}{\mathcal{Y}}
\nc{\J}{\mathcal{J}}
\nc{\cL}{\mathcal{L}}
\nc{\cM}{\mathcal{M}}
\nc{\cN}{\mathcal{N}}
\nc{\cH}{\mathcal{H}}
\nc{\cO}{\mathcal{O}}
\nc{\Op}{\mathrm{O}}
\nc{\rmA}{\mathrm{A}}
\nc{\rmF}{\mathrm{F}}
\nc{\rmB}{\mathrm{B}}
\nc{\rmD}{\mathrm{D}}
\nc{\rmP}{\mathrm{P}}
\nc{\rmT}{\mathrm{T}}
\nc{\cC}{\mathcal{C}}
\nc{\cP}{\mathcal{P}}
\nc{\bbQ}{\mathbb{Q}}
\nc{\bbR}{\mathbb{R}}
\nc{\cU}{\mathcal{U}}
\nc{\Un}{\bigcup}
\nc{\UnN}[1]{\Un_{#1=1}^\oo}
\nc{\CapN}[1]{\bigcap_{#1=1}^\infty}
\nc{\cV}{\mathcal{V}}
\nc{\cW}{\mathcal{W}}
\nc{\Z}{{\mathbb Z}}
\nc{\Impl}{\Rightarrow}
\long\def\forget#1\forgotten{\marginpar{\textcolor{green}{Forgetting...}}}
\nc{\ft}{\mathfrak{t}}
\nc{\fb}{\mathfrak{b}}
\nc{\fc}{\mathfrak{c}}
\nc{\fd}{\mathfrak{d}}
\nc{\fg}{\mathfrak{g}}
\nc{\oo}{\infty}
\nc{\fr}{\mathfrak{r}}
\nc{\fk}{\mathfrak{k}}
\nc{\bidi}{\mathfrak{bidi}}
\nc{\fu}{\mathfrak{u}}
\nc{\fh}{\mathfrak{h}}
\nc{\fp}{\mathfrak{p}}
\nc{\fj}{\mathfrak{j}}
\nc{\fs}{\mathfrak{s}}
\nc{\w}{\omega}
\nc{\x}{\times}
\nc{\Iff}{\Leftrightarrow}
\newcommand\comp{^{\text{\tt c}}}
\nc{\nin}{\notin}
\nc{\cat}{\hat{\ }}
\nc{\sub}{\subseteq}
\nc{\spst}{\supseteq}
\nc{\sm}{\setminus}
\nc{\as}{\subseteq^*}%{\let\proclaim\relax}
\nc{\les}{\le^*}
\nc{\leinf}{\le^{\infty}}
\nc{\leS}{\le_S}
\nc{\leF}{\le_{\mathcal{F}}}
\nc{\leU}{\le_{\mathcal{U}}}
\nc{\rest}{\restriction}

\nc{\la}{\langle}
\nc{\ra}{\rangle}
\nc{\E}{\exists}
\nc{\dom}{\op{dom}}
\nc{\cov}{\op{cov}}
\nc{\add}{\op{add}}
\nc{\cof}{\op{cof}}
\nc{\cf}{\op{cf}}
\nc{\non}{\op{non}}
\nc{\unif}{\op{non}}
\nc{\COV}{\op{COV}}
\nc{\ADD}{\op{ADD}}
\nc{\COF}{\op{COF}}
\nc{\NON}{\op{NON}}
\nc{\impl}{\to}
\nc{\Lp}{\mathcal{L_\p}}
\nc{\Wlog}{without loss of generality}
\newtheorem{thm}{Theorem}[section]
\nc{\bthm}{\begin{thm}} \nc{\ethm}{\end{thm}}
\newtheorem{prop}[thm]{Proposition}
\nc{\bprp}{\begin{prop}} \nc{\eprp}{\end{prop}}
\newtheorem{fact}[thm]{Fact}
\nc{\bfct}{\begin{fact}} \nc{\efct}{\end{fact}}
\newtheorem{prob}[thm]{Problem}
\nc{\bprb}{\begin{prob}} \nc{\eprb}{\end{prob}}
\newtheorem{lem}[thm]{Lemma}
\nc{\blem}{\begin{lem}} \nc{\elem}{\end{lem}}
\newtheorem{app}[thm]{Application}
\nc{\bapp}{\begin{app}} \nc{\eapp}{\end{app}}
\newtheorem{claim}[thm]{Claim}
\nc{\bclm}{\begin{claim}} \nc{\eclm}{\end{claim}}
\newtheorem{cor}[thm]{Corollary}
\nc{\bcor}{\begin{cor}} \nc{\ecor}{\end{cor}}
\newtheorem{conj}[thm]{Conjecture}
\nc{\bcnj}{\begin{conj}} \nc{\ecnj}{\end{conj}}
\theoremstyle{definition}
\newtheorem{defn}[thm]{Definition}
\nc{\bdfn}{\begin{defn}} \nc{\edfn}{\end{defn}}
\newtheorem{obs}[thm]{Observation}
\nc{\bobs}{\begin{obs}} \nc{\eobs}{\end{obs}}
\theoremstyle{remark}
\newtheorem{rem}[thm]{Remark}
\nc{\brem}{\begin{rem}} \nc{\erem}{\end{rem}}
\newtheorem{cnv}[thm]{Convention}
\nc{\bcnv}{\begin{cnv}} \nc{\ecnv}{\end{cnv}}
\newtheorem{exam}[thm]{Example}
\nc{\bexm}{\begin{exam}} \nc{\eexm}{\end{exam}}
\nc{\bpf}{\begin{proof}} \nc{\epf}{\end{proof}}
\nc{\be}{\begin{enumerate}}
\nc{\ee}{\end{enumerate}}
\nc{\bi}{\begin{itemize}}
\nc{\bimy}{\my{\begin{itemize}}
\nc{\eimy}{\end{itemize}}}
\nc{\itm}{\item}
\nc{\ei}{\end{itemize}}
\nc{\Subsection}[1]{\goodbreak\subsection*{#1}}%\ \par}

\nc{\sone}{\mathsf{S}_1}
\nc{\sfin}{\mathsf{S}_\mathrm{fin}}
\nc{\ufin}{\mathsf{U}_\mathrm{fin}}
\nc{\Split}{\mathsf{Split}}

\nc{\gone}{\mathsf{G}_1}    \nc{\gfin}{\mathsf{G}_\mathrm{fin}}

\newcommand{\ed}{

\end{document}
}

\title[Proofs from the book]{Selection principles and proofs from the Book}

\author{Boaz Tsaban}
\address{Boaz Tsaban,
Department of Mathematics,
Bar-Ilan University, Ramat Gan, Israel}
\email{tsaban@math.biu.ac.il}

\subjclass{%
Primary: 37F20; %Combinatorics and topology
Secondary 26A03, %Foundations: limits and generalizations, elementary topology of the line
03E75 %Applications of set theory
}

\begin{document}

\begin{abstract}
I provide simplified proofs for each of the following fundamental theorems regarding selection principles:
\be
\item The Quasinormal Convergence Theorem, due to the author and Zdomskyy,
asserting that a certain, important property of the space of continuous functions on a space is actually preserved by Borel images of that space.
\item The Scheepers Diagram Last Theorem, due to Peng,
completing all provable implications in the diagram.
\item The Menger Game Theorem, due to Telg\'arsky,
determining when Bob has a winning strategy in the game version of Menger's covering property.
\item A lower bound on the additivity of Rothberger's covering property, due to Carlson.
\ee
The simplified proofs lead to several new results.
\end{abstract}

\maketitle

\centerline{\emph{To Adina}}

\section{Introduction}

The study of \emph{selection principles} unifies
notions and studies originating from
dimension theory (Menger and Hurewicz), measure theory (Borel), convergence properties (Cs\'asz\'ar--Laczkovicz),
and function spaces (Gerlits--Nagy and Arhangel'ski\u{\i}),
notions analyzed and developed in numerous studies of later mathematicians, especially since the 1996 paper of Just, Miller, Scheepers and Szeptycki~\cite{coc2}.
The unified notions include, among others, many classic types of special sets of
real numbers, local properties in function spaces,
and more recent types of convergence properties.

\begin{figure}[!htp]
	\begin{changemargin}{-4cm}{-3cm}
		\begin{center}
			{%\scriptsize
				$\xymatrix@R=8pt{
					%1
					&
					&
					& \sr{$\ufin(\Op,\Ga)$}{Hurewicz}{$\fb$}\ar[r]
					& \sr{$\ufin(\Op,\Om)$}{}{$\fd$}\ar[rr]
					& & \sr{$\sfin(\Op,\Op)$}{Menger}{$\fd$}
					\\
					%2
					&
					&
					& \sr{$\sfin(\Ga,\Om)$}{}{$\fd$}\ar[ur]
					\\
					%3
					& \sr{$\sone(\Ga,\Ga)$}{}{$\fb$}\ar[r]\ar[uurr]
					& \sr{$\sone(\Ga,\Om)$}{}{$\fd$}\ar[rr]\ar[ur]
					& & \sr{$\sone(\Ga,\Op)$}{}{$\fd$}\ar[uurr]
					\\
					%4
					&
					&
					& \sr{$\sfin(\Om,\Om)$}{}{$\fd$}\ar'[u][uu]
					\\
					%5
					& \sr{$\sone(\Om,\Ga)$}{Gerlits--Nagy}{$\fp$}\ar[r]\ar[uu]
					& \sr{$\sone(\Om,\Om)$}{}{$\cov(\cM)$}\ar[uu]\ar[rr]\ar[ur]
					& & \sr{$\sone(\Op,\Op)$}{Rothberger}{$\cov(\cM)$}\ar[uu]
				}$
			}
			\caption{The Scheepers Diagram}\label{SchDiag}
		\end{center}
	\end{changemargin}
\end{figure}
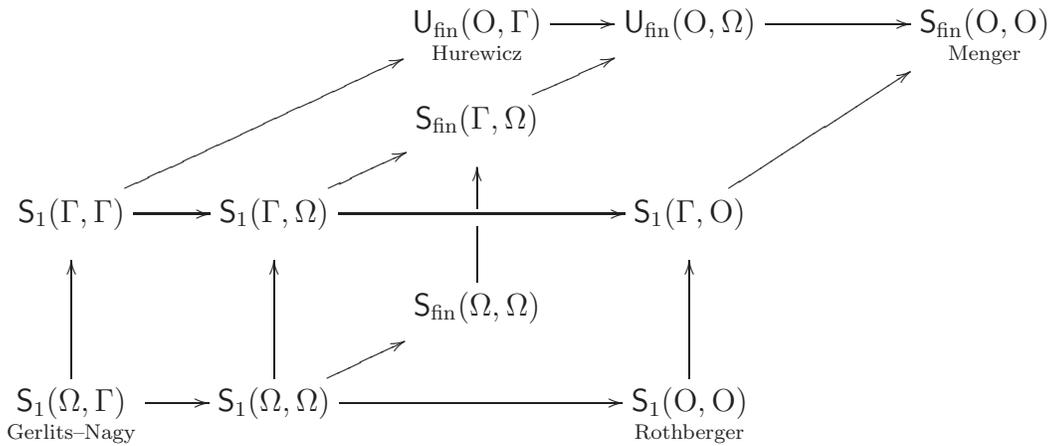

Selective topological covering properties form the kernel of selection principles.
These covering properties are related via the Scheepers Diagram (Figure~\ref{SchDiag}).
This is a diagram of covering properties and implications among them.
The properties in this diagram are obtained as follows.

For families $\rmA$ and $\rmB$ of sets,
let $\sone(\rmA,\rmB)$ be the statement:
For each sequence of elements of the family $\rmA$, we can pick one element from each sequence
member, and obtain an element of the family $\rmB$.
When $\rmA=\rmB=\Op(X)$, the family of open covers of a topological space $X$, we obtain
\emph{Rothberger's property} (1941), the topological version of Borel's strong measure zero.
We say that a space $X$ satisfies $\sone(\Op,\Op)$ if the assertion $\sone(\Op(X),\Op(X))$ holds,
and similarly for the other selective properties.

The hypothesis $\sfin(\rmA,\rmB)$ is obtained by replacing \emph{one} by
\emph{finitely many} in the above definition. The property $\sfin(\Op,\Op)$ is, by an observation
of Hurewicz (1925), equivalent to Menger's basis property, a dimension-type property.
The property $\ufin(\rmA,\rmB)$ is obtained by further allowing us to take the \emph{unions}
of the selected finite subsets---this matters for some types of covers.
For technical reasons, this property does not consider all covers of type $\rmA$, but only those
that have no finite subcover.

A cover of a space is an \emph{$\omega$-cover} if no member of the cover covers the entire space, but
every finite subset of the space is covered by some member of the cover.
For a space $X$, $\Om(X)$ is the family of open $\omega$-covers of the space.
A  \emph{point-cofinite cover} is an infinite cover where every point of the space belongs to all but
finitely many members of the cover.
$\Ga(X)$ is the family of open point-cofinite covers of the space $X$.

Applying the mentioned selection principles to the cover types $\Op$, $\Om$ and $\Ga$,
we obtain additional important properties, such as Hurewicz's property $\ufin(\Op,\Ga)$ (1925).
We also obtain the Gerlits--Nagy $\gamma$-property $\sone(\Om,\Ga)$ (1982), characterizing
the Fr\'echet--Urysohn property of the function space $\Cp(X)$ of continuous real-valued functions,
with the topology of pointwise convergence:
A topological space is \emph{Fr\'echet--Urysohn} if every point in the closure of
a set is actually a limit of a sequence in the set.
This duality between the spaces $X$ and $\Cp(X)$ also translates
various tightness and convergence properties of the space $\Cp(X)$---discovered
earlier by Arhangel'ski\u{\i}, Bukovsk\'y, Sakai, and others---to the selective covering properties
$\sfin(\Om,\Om)$, $\sone(\Ga,\Ga)$, and $\sone(\Om,\Om)$.
In Section~\ref{sec:QN}, we provide a surprisingly simple proof of one of the most important results of this type.
While the result itself does not involve selective covering properties explicitly,
its proof does that extensively.

A topological space is \emph{Lindel\"of} if every open cover has a countable subcover.
For example, all sets of real numbers are Lindel\"of.
Since all selection principles concern countable sequences, the theory mainly deals with
Lindel\"of spaces.
For Lindel\"of spaces, the Scheepers Diagram is the result of a classification of all properties thus introduced;
each property is equivalent to one in the diagram~\cite{coc2}.
It was long open whether any additional implication could be established among the properties in the diagram.
In Section~\ref{sec:Peng} we deal with the recent, surprising solution of this problem.

Menger's covering property $\sfin(\Op,\Op)$ is the oldest, most general, and most applied property in the Scheepers Diagram.
Initially, Menger conjectured his property to coincide
with $\sigma$-compactness.
While this turned out false~\cite{MHP}, the game version of this
property does provide a characterization of $\sigma$-compactness.
A very transparent proof of this deep result is presented in Section~\ref{sec:Bob}.

In Section~\ref{sec:add} we consider a connection to combinatorial set theory.
We show that a nontrivial lower bound on the additivity of Rothberger's property
follows easily from basic knowledge on selection principles.

\emph{The Book} is a popular myth by Paul Erd\H{o}s: A transfinite book containing the most simple proofs for all theorems.
I would like to believe that the proofs presented here are similar to ones from the Book \dots or from some of its preliminary drafts, at any rate.
%Perhaps humankind is taking part in the writing of the Book;
%I could not access  Erd\H{o}s to verify any of these hypotheses.

\section{The Quasinormal Convergence Theorem}
\label{sec:QN}

By \emph{real set} we mean a topological space where every open set is a countable union of clopen sets.
Such are, for example, totally disconnected subsets of the real line and, in particular, subsets of the Cantor space $\Cantor$.
In general, every perfectly normal space with any of the properties considered in this section is a real set.

Let $X$ be a real set.
A sequence of real-valued functions $\eseq{f}$ on $X$ converges \emph{quasinormally}
to a real-valued function $f$
if there are positive real numbers $\eseq{\epsilon}$ converging to $0$
such that for each point $x\in X$, we have
\[
\card{f_n(x)-f(x)}\le \epsilon_n
\]
for all but finitely many $n$.
Quasinormal convergence generalizes uniform convergence.

A real set $X$ is a \emph{QN space} if
every sequence of continuous real-valued functions on $X$ that converges to $0$ pointwise,
converges to $0$ quasinormally. Equivalently, convergence in the space $\Cp(X)$
is quasinormal.

QN spaces were studied intensively, e.g., by Bukovsk\'y, Rec\l{}aw, Repick\'y, Scheepers, Nowik, Sakai, and Hale\v{s}~\cite[and references therein]{hH}.
This and other properties of similar type are preserved by continuous images, and all experience prior to the
paper of the author and Zdomskyy~\cite{hH}
suggested that they are not preserved by Borel images.
Thus, the following theorem~\cite[Theorem~9]{hH} came as a surprise.

The Baire space $\NN$ is quasiordered by the relation $\les$ of eventual dominance:
$f\les g$ if $f(n)\le g(n)$ for all but finitely many $n$.
A subset $Y$ of the Baire space $\NN$ is \emph{bounded} if it is bounded with respect to evntual dominance.

\bthm[Quasinormal Convergence Theorem]
The following assertions are equivalent for real sets $X$:
\be
\itm The set $X$ is a QN space.
\itm Every Borel image of the set $X$ in the Baire space $\NN$ is bounded.
\ee
\ethm

The second property in the theorem is well known and straightforward to apply:
The most natural transformations needed in proofs regarding these notions are always easily seen to be Borel.
Consequently, the theorem had a dramatic impact on the study of QN spaces:
First, many of the previous sophisticated arguments could be replaced by
straightforward ones.
Second, many properties that were hitherto considered separately turned out provably equivalent.
Consequently, this theorem settled all problems concerning these properties~\cite{hH}.

The original proof of the Quasinormal Convergence Theorem is long and involved, and some of its parts are difficult to follow.
A more natural proof was later published by Bukovsk\'y and \v{S}upina~\cite[\S4]{BukSup12}.
Inspired by a paper of Gerlits and Nagy~\cite{GNFr}, I have discovered the following surprisingly simple proof.
All needed proof ingredients were already available at the time the Quasinormal Convergence Theorem was established.
The following lemma provides the key to the proof.

For a space $X$, let $\PF(X)$ be the collection of countably infinite point-finite families of open sets in $X$.

\blem
\label{lem:PF}
Let $X$ be a topological space. The following assertions are equivalent:
\be
\item Every Borel image of the space $X$ in the Baire space $\NN$ is bounded.
\item The space $X$ satisfies $\sone(\PF,\PF)$.
\ee
\elem
\bpf
Let $\rmF(X)$ (respectively, $\rmB(X)$) be the family of countable closed (respectively, Borel) covers of the set $X$,
and $\FG(X)$ (respectively, $\BG(X)$) be the family of infinite closed (respectively, Borel) point-cofinite covers of the set $X$.
The properties (1), $\ufin(\rmB,\BG)$, and $\sone(\BG,\BG)$ are equivalent~\cite[Theorem~1]{cbc}.
For a family $\cU$ of open sets, we have $\cU\in\PF(X)$ if and only if
\[
\sset{U\comp}{U\in\cU}\in\FG(X).
\]
It follows that $\sone(\PF,\PF)=\sone(\FG,\FG)$.

$(1)\Impl (2)$: Clearly, $\sone(\BG,\BG)$ implies $\sone(\FG,\FG)$.

$(2)\Impl (1)$: A theorem of Bukovsk\'y--Rec{\l}aw--Repick\'y~\cite[Corollary~5.3]{BRR91} asserts that
\[
\ufin(\rmF,\FG)=\ufin(\B,\BG).
\]
The usual argument~\cite[Proposition~11]{coc1} shows that $\sone(\FG,\FG)$ implies $\ufin(\rmF,\FG)$:
If $\sseq{C_n}\in\rmF(X)$ and there is no finite subcover, then $\sseq{\Un_{k=1}^n C_k}\in\FG(X)$.
\epf

A topological space $Y$ has Arhangel'ski\u{\i}'s property $\alpha_1$ if
for every sequence $\eseq{s}$ of sequences converging to the same point,
there is a sequence $s$ such that
the sets $\im(s_n)\sm \im(s)$ are finite for all natural numbers $n$.
This property is defined by properties of sets (images of sequences) rather than sequences.
Fix a bijection $\varphi\colon \bbN\x\bbN\to\bbN$.
For sequences $\eseq{s}$, with $s_n=(s_{(n,1)},s_{(n,2)},\dotsc)$ for each $n$,
define
\[
\bigsqcup_{n=1}^\oo s_n:=(s_{\varphi(1)},s_{\varphi(2)},\dotsc).
\]
Since convergence of a sequence does not depend on the order of its elements, it does not matter,
for our purposes, which bijection $\varphi$ is used.
A sequence
$\bigsqcup_{n=1}^\oo s_n$
converges to a point $p$ if and only if each sequence $s_n$ converges to $p$,
and for each neighborhood $U$ of $p$, we have $\im(s_n)\sub U$ for all but finitely many $n$.

\blem
Let $Y$ be an $\alpha_1$ space.
For every sequence $\eseq{s}$ of sequences in the space $Y$ converging to the same point $p$,
there are tails $t_n$ of $s_n$, for $n\in\bbN$, such that the sequence
$\bigsqcup_{n=1}^\oo t_n$ converges to $p$.
\elem
\bpf
There is a sequence $s$ such that
the sets $\im(s_n)\sm \im(s)$ are finite for all natural numbers $n$.
By moving to a subsequence, we may assume that
$\im(s)\sub\UnN{n}\im(s_n)$.
Suppose that $s=(\eseq{a})$.
For each natural number $n$, since the sequence $s_n$ coverges to the point $p$,
every element other than $p$ may appear in the sequence $s_n$ only finitely often.
Thus, there is a tail $t_n$ of the sequence $s_n$ such that
\[
\im(t_n)\sub\sset{a_k}{k\geq n}\cup\{p\}.
\]
Let $U$ be a neighborhood of $p$.
There is a natural number $N$ such that
\[
\im(t_n)\sub\sset{a_k}{k\geq n}\cup\{p\}\sub U
\]
for all natural numbers $n\ge N$.
Thus, the direct sum $t:=\bigsqcup_{n=1}^\oo t_n$ converges to the point $p$.
\epf

Sakai~\cite[Theorem~3.7]{Sakai07} and Bukovsk\'y--Hale\v{s}~\cite[Theorem~11]{BH07} proved
that a real set $X$ is a QN space if, and only if, the space $\Cp(X)$ is an $\alpha_1$ space.
Thus, the Quasinormal Convergence Theorem can be stated, and proved, as follows.

\bthm
\label{thm:pf}
The following assertions are equivalent for real sets $X$:
\be
\itm The space $C_p(X)$ is an $\alpha_1$ space.
\itm Every Borel image of the set $X$ in the Baire space $\NN$ is bounded.
\ee
\ethm
\bpf
$(2)\Impl (1)$: This is the straightforward implication.
For completeness, we reproduce its proof~\cite[Theorem~9]{hH}.

Let $\eseq{s}$ be sequences in the space $\Cp(X)$ that converge to a function $f\in \Cp(X)$.
For each natural number $n$, suppose that
\[
s_n = (f^n_1,f^n_2,f^n_3,\dotsc).
\]
Define a Borel function $\Psi\colon X\to\NN$ by
\[
\Psi(x)(n) := \min\medset{k}{(\forall m\ge k)\ |f^n_m(x)-f(x)|\le \frac{1}{n}}.
\]
Let $g\in\NN$ be a $\les$-bound for the image $\Psi[X]$.
Then the sequence
\[
\bigsqcup_{n=1}^\oo(f^n_{g(n)},f^n_{g(n)+1},f^n_{g(n)+2},\dotsc)
\]
converges to the function $f$.

$(1)\Impl (2)$: This is the main implication. By Lemma~\ref{lem:PF}, it suffices to prove that the
set $X$ satisfies $\sone(\PF,\PF)$.
Let $\eseq{\cU}\in\PF(X)$.
By thinning out the point-finite covers, we may assume that they are pairwise disjoint~\cite[Lemma~4]{coc1}.
For each set $U\in\UnN{n} \cU_n$, let $\cC_U$ be a countable family of clopen sets with $\Un\cC_U=U$.
%By the standard diagonal argument~\cite[Lemma~4]{coc1}, taking finite unions inside each family $\cC_U$,
%we may assume that the families $\cC_U$ are pairwise disjoint.
For each natural number $n$, let
\[
\cV_n := \Un_{U\in\cU_n}\cC_U.
\]
Every set $C\in\cV_n$ is contained in at most finitely many sets $U\in\cU_n$.
Thus, the family $\cV_n$ is infinite and point-finite.
Let $s_n$ be a bijective enumeration of the family
\[
\sset{\chi_V}{V\in\cV_n}.
\]
The sequence $s_n$ is in $\Cp(X)$, and it converges to the constant function $0$.

As the space $\Cp(X)$ is $\alpha_1$, there are for each $n$ a tail $t_n$ of the sequence $s_n$
such that the sequence $s:=\bigsqcup_{n=1}^\oo t_n$ converges to $0$.
%It follows that the family $\cC:=\sset{V}{\chi_V\in \im(s)}$ is point-finite.
For each natural number $n$, pick a set $U_n\in \cU_n$ with $\cC_{U_n}\sub\im(t_n)$.
Then the family $\sseq{U_n}$ is infinite and point-finite.
\epf

For a set $X\sub\Cantor$, Gerlits and Nagy (and, independently, Nyiko\v{s}) define a space $\rmT(X)$ as follows.
Let $\{0,1\}^*$ denote the set of finite sequences of elements of the set
$\{0,1\}$.
Let $X\sub\Cantor$.
For each point $x\in X$, let $A_x\sub\{0,1\}^*$ be the
set of initial segments of the point $x$.
Let $X\cup\{0,1\}^*$
be the topological space where the points of the set $\{0,1\}^*$ are isolated,
and for each point $x\in X$, a neighborhood base of $x$ is given by the sets $\{x\}\cup B$,
where $B$ is a cofinite subset of the set $A_x$.
Let $\rmT(X)$ be the one-point compactification of this space, and $\oo$ be the compactifying point.

Gerlits and Nagy prove that if a set $X\sub\Cantor$ is a Siepi\'nski set, then
the space $\rmT(X)$ is $\alpha_1$, and that if the space $\rmT(X)$ is $\alpha_1$, then the set $X$ is a $\sigma$-set~\cite[Theorem~4]{GNFr} .
%In the concluding Theorem of their paper, they prove the same assertions when the space $\rmT(X)$ is replaced by the space $\Cp(X)$.
The following theorem unifies these results and improves upon them.
Indeed, every Borel image of a Siepi\'nski set in the Baire space is bounded,
and every set with bounded Borel images in the Baire space  is a $\sigma$-set~\cite[and references therein]{cbc}.

\bthm
Let $X\sub\Cantor$. The following assertions are equivalent:
\be
\itm The space $\rmT(X)$ is an $\alpha_1$ space.
\itm Every Borel image of the set $X$ in the Baire space $\NN$ is bounded.
\ee
\ethm
\bpf
For  a finite sequence $s\in\{0,1\}^*$, let $[s]$ be the basic clopen subset of the Cantor space $\Cantor$ consisting of all functions
extending $s$.
Every open set in the space $\Cantor$ is a disjoint union of basic clopen sets.
A sequence $\eseq{a}$ in the set $\{0,1\}^*$ converges to $\oo$ in the space $\rmT(X)$ if, and only if,
the set $\sseq{[a_n]}$ is point-finite in the space $X$~\cite{GNFr}.
The argument in the proof of Theorem~\ref{thm:pf} applies.
\epf

\section{The Scheepers Diagram Last Theorem}
\label{sec:Peng}

\begin{figure}[!htp]
	\begin{changemargin}{-4cm}{-3cm}
		\begin{center}
			{%\scriptsize
				$\xymatrix@R=8pt{
					%2
					&
					& \sr{$\ufin(\Op,\Ga)$}{Hurewicz}{$\fb$}\ar[r]
					& \sr{$\sfin(\Ga,\Om)$}{}{$\fd$}\ar[r]
					& \sr{$\ufin(\Op,\Om)$}{}{$\fd$}\ar[r]
					& \sr{$\sfin(\Op,\Op)$}{Menger}{$\fd$}
					\\
					%3
					& \sr{$\sone(\Ga,\Ga)$}{}{$\fb$}\ar[r]\ar[ur]
					& \sr{$\sone(\Ga,\Om)$}{}{$\fd$}\ar[rr]\ar[ur]
					& & \sr{$\sone(\Ga,\Op)$}{}{$\fd$}\ar[ur]
					\\
					%4
					&
					&
					& \sr{$\sfin(\Om,\Om)$}{}{$\fd$}\ar'[u][uu]
					\\
					%5
					& \sr{$\sone(\Om,\Ga)$}{Gerlits--Nagy}{$\fp$}\ar[r]\ar[uu]
					& \sr{$\sone(\Om,\Om)$}{}{$\cov(\cM)$}\ar[uu]\ar[rr]\ar[ur]
					& & \sr{$\sone(\Op,\Op)$}{Rothberger}{$\cov(\cM)$}\ar[uu]
				}$
			}
			\caption{The Final Scheepers Diagram}\label{fig:SchDiagFinal}
		\end{center}
	\end{changemargin}
\end{figure}
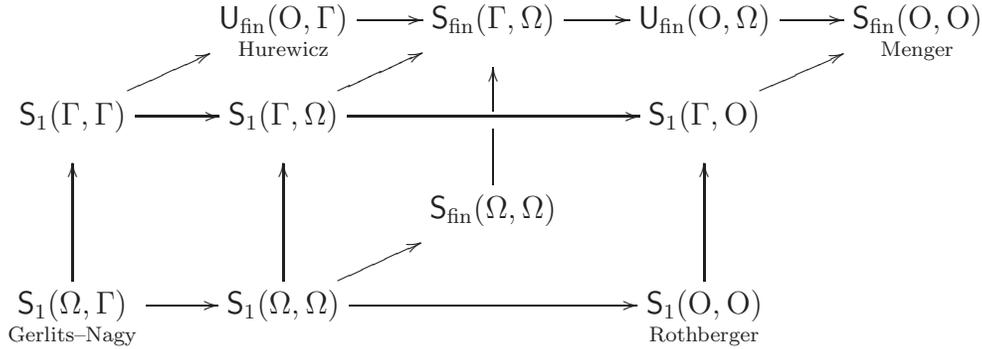

The implications in the Scheepers Diagram~\ref{SchDiag} were all rather straightforward to establish, and
almost all other potential implications were ruled out by counterexamples~\cite{coc2}.
Only two problems remained open: Does $\ufin(\Op,\Om)$ imply $\sfin(\Ga,\Om)$?
And if not, does $\ufin(\Op,\Ga)$ imply $\sfin(\Ga,\Om)$?~\cite[Problems~1 and~2]{coc2}.
For nearly three decades it was expected that the remaining two potential implications were refutable.
Only when Peng came up with an entirely new method for refuting implications among selective covering properties~\cite{Peng},
were these problems resolved. But not in the expected way: Having proved that $\ufin(\Op,\Om)$ does not imply $\sfin(\Ga,\Om)$,
Peng tried to refute the last remaining potential implication. And he failed.
His close examination of the failure suggested a path for \emph{proving} the last potential implication~\cite[Theorem~23]{Peng}.
Peng's results establish the final form of the Scheepers Diagram (Figure~\ref{fig:SchDiagFinal}).

Peng's proof of the last implication is somewhat involved.
The proof given below identifies the heart of Peng's argument, and replaces
the other parts with simple, quotable observations about selective covering properties.

Let $k$ be a natural number.
A cover of a space is a \emph{$k$-cover} if no member of the cover covers the entire space, but
every $k$-element subset of the space is covered by some member of the cover.
Thus, a cover is an $\w$-cover if and only if it is a $k$-cover for all natural numbers $k$.
For a space $X$ and a natural number $k$, let $\Op_k(X)$  be the family of open $k$-covers of the space.

\blem
\label{lem:ncov}
Let $\Pi$ be a selection principle, and $\rmA$ a type of open covers.
The following assertions are equivalent:
\be
\item $\Pi(\rmA,\Om)$.
\item For each natural number $k$ we have $\Pi(\rmA,\Op_k)$.
\ee
\elem
\bpf
$(1)\Impl(2)$: Obvious.

$(2)\Impl (1)$: Let $\eseq{\cU}$ be a sequence in $\rmA$. Split the sequence into infinitely many disjoint sequences.
For each natural number $k$, apply $\Pi(\rmA,\Op_k)$ to the $k$th sequence, to obtain a $k$-cover $\cV_k$. Then
$\UnN{k} \cV_k$ is an $\omega$-cover, in accordance to the required property $\Pi(\rmA,\Om)$.
\epf

\bthm[{Peng~\cite[Theorem~23]{Peng}}]
\label{thm:peng}
The Hurewicz property $\ufin(\Op,\Ga)$ implies $\sfin(\Ga,\Om)$.
\ethm
\bpf
Let $X$ be a Huewicz space.
By Lemma~\ref{lem:ncov}, it suffices to prove that $\sfin(\Ga,\Op_k)$ holds for all natural numbers $k$.
Fix a natural number $k$.

Let $\eseq{\cU}$ be a sequence in $\Ga(X)$.
By moving to countably infinite subcovers, we may enumerate
\[
\cU_n = \sset{U^n_m}{m\in\bbN}
\]
for each $n$.
For each $n$ and $m$, we may replace the set $U^n_m$ with the smaller set
\[
U^1_m\cap U^2_m\cap\dotsb\cap U^n_m,
\]
so that we may assume that
\[
U^1_m\spst U^2_m\spst U^3_m\spst\dotsb
\]
for all natural numbers $m$.
The refined covers $\cU_n$ remain in $\Ga(X)$.

Let $g_0(m):=m$ for all $m$.
We define, by induction, increasing functions $g_1,\dots,g_k\in\NN$.
Let $l<k$ and assume that the function $g_l$ is defined.
For natural numbers $n$, $m$ and $i$, let
\beq
V^{l,n}_i &:=& \bigcap_{m=i}^{g_l(i)}U^n_m;\\
W^{l,n}_m &:=& \Un\sset{V^{l,n}_i}{n\le i, g_l(i)\le m}.
\eeq
For each $l$ and $n$, the sets $W^{l,n}_m$ are increasing with $m$, and cover the space $X$.
By the Hurewicz property, there is an increasing function $g_{l+1}\in\NN$ such that
\[
\sseq{W^{l,n}_{g_{l+1}(n)}}\in\Ga(X).
\]
This completes the iductive construction.

We will show that
\[
\sset{U^n_m}{n\in\bbN, m\le g_k(n)}\in\Op_k(X).
\]
Let $x_1,\dotsc,x_k\in X$.
Since $\sseq{W^{l,n}_{g_{l+1}(n)}}\in\Ga(X)$ for all $l=0,\dotsc,k-1$,
there is a natural number $N$ with
\[
x_1,\dots,x_k\in W^{l,n}_{g_{l+1}(n)}
\]
for all $l=1,\dotsc,k$ and all $n\ge N$.
Fix a number $n_0\ge N$.

Since
$x_1\in W^{k-1,n_0}_{g_{k}(n_0)}$,
there is $n_1$ with
$n_0\le n_1,g_{k-1}(n_1)\le g_{k}(n_0)$ and
\[
x_1\in V^{k-1,n_0}_{n_1}=\bigcap_{m=n_1}^{g_{k-1}(n_1)}U^{n_0}_m.
\]
\nc{\sentence}[2]{
Since
$x_#2\in W^{k-#2,n_#1}_{g_{k-#1}(n_#1)}$,
there is $n_#2$ with
$n_#1\le n_#2,g_{k-#2}(n_#2)\le g_{k-#1}(n_#1)$ and
\[
x_#2\in V^{k-#2,n_#1}_{n_#2}=\bigcap_{m=n_#2}^{g_{k-#2}(n_#2)}U^{n_#1}_m\sub \bigcap_{m=n_#2}^{g_{k-#2}(n_#2)}U^{n_0}_m.
\]
}
\sentence{1}{2}
\sentence{2}{3}
\centerline{\vdots}

Since
$x_k\in W^{0,n_{k-1}}_{g_{1}(n_{k-1})}$,
there is $n_k$ with
$n_{k-1}\le n_k=g_{0}(n_k)\le g_{1}(n_{k-1})$ and
\[
x_k\in V^{0,n_{k-1}}_{n_k}=U^{n_{k-1}}_{n_k}\sub U^{n_0}_{n_k}.
\]
It follows that $x_1,\dots,x_k\in U^{n_0}_{n_k}$, and $n_k\le g_k(n_0)$.
\epf

The proof of Theorem~\ref{thm:peng} establishes a stronger result. To this end, we need the following definitions and lemma.
An infinite cover of a space $X$ is \emph{$\omega$-groupable} (respectively, \emph{$k$-groupable}, for a natural number $k$) if
there is a partition of the cover into finite parts such that
for each finite (respectively, $k$-element) set $F\sub X$  and all but finitely many parts $\cP$ of the partition, there is a set $U\in\cP$ with $F\sub U$~\cite{coc7}.
Let $\Om^\grbl(X)$ (respectively, $\Op_k^\grbl(X)$) be the family of open $\omega$-groupable (respectively, $k$-groupable) covers of the space $X$.

\blem
\label{lem:kcov2}
Let $\Pi$ be a selection principle, and $\rmA$ a type of open covers.
The following assertions are equivalent:
\be
\item $\Pi(\rmA,\Om^\grbl)$;
\item For each natural number $k$ we have $\Pi(\rmA,\Op_k^\grbl)$.
\ee
\elem
\bpf
The proof is similar to that of Lemma~\ref{lem:ncov}, once we observe that if $\sseq{U_n}$ is a $k$-grouable cover for all $k$, then
it is $\omega$-groupable.
This follows easily from the fact that for each countable family $\sset{\cP_k}{k\in\bbN}$ of partitions of $\bbN$ into finite sets,
there is a partition $\cP$ of $\bbN$ into finite sets that is eventually coarser than all of the given partitions, that is,
such that for each $k$, all but finitely many members of the partition $\cP$ contain a member of the partition $\cP_k$.
\epf

Ko\v{c}inac and Scheepers proved that if all finite powers of a space $X$ are Hurewicz, then every open $\omega$-cover of the space
is $\omega$-groupable. Together with Peng's Theorem~\ref{thm:peng}, we have that if \emph{all finite powers} of a space $X$ are Hurewicz,
then the space satisfies $\sfin(\Ga,\Om^\grbl)$. The following theorem shows that the assumption on the finite powers is not needed.

In the theorem, we also mention $\sfin(\Ga,\Lambda^\grbl)$.
An open cover is \emph{large} if each point is in infinitely many members of the cover.
Let $\Lambda(X)$ be the family of large open covers of the space $X$.
An open cover $\cU$ is in $\Lambda^\grbl(X)$~\cite{coc7} (also denoted
$\gimel(\Ga)$, depending on the context~\cite{GlCovs}) if there is a partition of the cover into finite parts such that
for each point $x\in X$ and all but finitely many parts $\cP$ of the partition, we have $x\in\Un\cP$.

\bthm
The following assertions are equivalent:
\be
\item $\ufin(\Op,\Ga)$,
\item $\sfin(\Ga,\Om^\grbl)$; and
\item $\sfin(\Ga,\Lambda^\grbl)$.
\ee
\ethm
\bpf
$(1)\Impl (2)$:
The proof of Peng's Theorem~\ref{thm:peng}, as written above,
shows, for a prescribed number $k$, that for each $k$-element set $F$ there is a natural number $N$ such that for each $n\ge N$ there is
a member of the finite set $\cF_{n}=\sset{U^{n}_m}{m\le g_k(n)}$ that contains the set $F$.
By thinning out the point-cofinite covers, we may assume that they are pairwise disjoint~\cite[Lemma~4]{coc1},
and consequently so are the finite sets $\cF_{n}$. Thus, $\UnN{n}\cF_n\in \Op_k^\grbl$.
This proves $\sfin(\Ga,\Op_k^\grbl)$ for all $k$.
Apply Lemma~\ref{lem:kcov2}.

$(2)\Impl (3)$: $\Om^\grbl\sub \Lambda^\grbl$.

$(3)\Impl (1)$:
This implication is standard and should be known. For completeness, we provide a proof.
Assume that the space $X$ satisfies $\sfin(\Ga,\Lambda^\grbl)$.
It suffices to prove that it satisfies $\ufin(\Ga,\Ga)$.
Given a sequence $\eseq{\cU}$ in $\Ga(X)$, we may (as in the proof of Theorem~\ref{thm:peng})
assume that the covers get finer with $n$.
Apply $\sfin(\Ga,\Lambda^\grbl)$ to obtain a cover $\cU\in \Lambda^\grbl$, with parts $\cP_n$ (for $n\in\bbN$) witnessing that.

Let $\cF_1\sub\cU_1$ be a finite set refined by $\cP_1$.
Let $n_2$ be minimal with $\cP_{n_2}\sub \Un_{n=2}^\oo \cU_n$, and
$\cF_2\sub\cU_2$ be a finite set refined by $\cP_{n_2}$.
Let $n_3$ be minimal with $\cP_{n_3}\sub \Un_{n=3}^\oo \cU_n$,
and $\cF_3\sub\cU_3$ be a finite set refined by $\cP_{n_3}$.
Continuing in this manner, we obtain finite sets $\cF_n\sub\cU_n$ for $n\in\bbN$, with
$\sseq{\Un\cF_n}\in\Ga(X)$.
\epf

Ko\v{c}inac and Scheepers proved that $\ufin(\Op,\Ga)=\sfin(\Om,\Lambda^\grbl)=\sfin(\Lambda,\Lambda^\grbl)$~\cite[Theorem~14]{coc7}.
However, $\ufin(\Op,\Ga)\neq \sfin(\Om,\Om^\grbl)$: The latter property is equivalent to satisfying $\ufin(\Op,\Ga)$ in all
finite powers~\cite[Theorem~16]{coc7}, a property strictly stronger than $\ufin(\Op,\Ga)$~\cite[Theorem~2.12]{coc2}.

\section{When Bob has a winning strategy in the Menger game}
\label{sec:Bob}

Menger~\cite{Menger24} conjectured that his property $\sfin(\Op,\Op)$ implies $\sigma$-compactness.
While his conjecture turned out false~\cite[and references therein]{MHP}, a closely related assertion is true.
The \emph{Menger game}~\cite{Hure25}, $\gfin(\Op,\Op)$, is the game associated to Menger's property $\sfin(\Op,\Op)$.
It is played on a topological space $X$, and has an inning per each natural number $n$.
In each inning, Alice picks an open cover $\cU_n$ of the space, and Bob chooses a finite set $\cF_n\sub\cU_n$.
Bob wins if $\UnN{n}\cF_n$ is a cover of the space, and otherwise Alice wins.
Telg\'arsky~\cite{TelGames3} proved that if Bob has a winning strategy in the Menger game played on a metric space,
then the space  is $\sigma$-compact.

Scheepers~\cite[Theorem~1]{Sch95} provided a direct proof of Telg\'arsky's Theorem, using the notion of H-closed sets.
We will eliminate the notion of H-closed sets and the closure operations from Scheepers's proof, and obtain a more transparent proof.
As a bonus, the separation hypotheses on the space are eliminated.

A subset $K$ of a topological space $X$ is \emph{relatively compact} if
every open cover $\cU$ of the entire space $X$ has a finite subcover of the set $K$.
A set $K$ is relatively compact if and only if its closure is compact.

\blem
\label{lem:kappacompact}
Let $\kappa$ be a cardinal number.
If a space $X$ is a union of at most $\kappa$ relatively compact sets, then it is the union of at most $\kappa$ compact sets.
\elem
\bpf
If $X=\Un_{\alpha<\kappa} K_\alpha$, then $X=\Un_{\alpha<\kappa} \cl{K_\alpha}$.
\epf

For a basis $\cB$ for the topology of a space $X$, let $\Op_\cB(X)$ be the family of subsets of $\cB$ that cover the space $X$.

\blem
\label{lem:relcompact}
Let $X$ be a topological space with a basis $\cB$, and $\sigma$ be a function on the family $\Op_\cB(X)$ such that for each cover $\cU\in\Op_\cB(X)$,
$\sigma(\cU)$ is a finite subset of $\cU$. Then the set
\[
K:=\bigcap_{\cU\in\Op_\cB(X)}\Un\sigma(\cU)
\]
is relatively compact.
\elem
\bpf
Let $\cU$ be an open cover of $X$. Let $\cV\in\Op_\cB(X)$ be a cover that refines the cover $\cU$.
Then $K\sub \Un\sigma(\cV)$, and there is a finite set $\cF\sub\cU$ with $ \Un\sigma(\cV)\sub\Un\cF$.
\epf

Scheepers~\cite[Theorem~1]{Sch95} proves the following theorem for metric spaces.
If Bob has a winning strategy in the Menger game played on $X$, then the space $X$ is Menger and, in particular, Lindel\"of.
If $X$ is, in addition, metric, then the space is second countable.

\bthm
\label{thm:Mgame1}
Let $X$ be a second countable topological space. If Bob has a winning strategy in the Menger game
$\gfin(\Op,\Op)$
played on $X$, then the space $X$ is $\sigma$-compact.
\ethm
\bpf
We follow steps of Scheepers's proof, removing what is not necessary.
Let $\sigma$ be a winning stratery for Bob.
Fix a countable base $\cB$ for the topology of the space $X$.
Let $\bbN^*$ be the set of finite sequences of natural numbers.
We consider all possible games where Alice chooses her covers from the family $\Op_\cB(X)$.

Since the base $\cB$ is countable, the family $\sset{\sigma(\cU)}{\cU\in\Op_\cB}$ (the possible first responds of Bob) is countable, too.
Choose elements $\eseq{\cU}\in\Op_\cB$ with
\[
\sseq{\sigma(\cU_n)}=\sset{\sigma(\cU)}{\cU\in\Op_\cB}.
\]
By induction, for a given natural number $n$ and each sequence $s\in\bbN^n$, the family
\[
\sset{\sigma(\cU_{s_1},\cU_{s_1,s_2},\dotsc,\cU_s,\cU)}{\cU\in\Op_\cB}
\]
is countable. Choose elements $\cU_{s,1},\cU_{s,2},\dotsc\in\Op_\cB$ with
\[
\sseq{\sigma(\cU_{s_1},\dotsc,\cU_{s},\cU_{s,n})}
=
\sset{\sigma(\cU_{s_1},\dotsc,\cU_{s},\cU)}{\cU\in\Op_\cB}.
\]
This completes our inductive construction.

By Lemma~\ref{lem:relcompact}, for each sequence $s\in\bbN^*$,  the set
\[
K_s := \CapN{n}\Un\sigma(\cU_{s_1},\dotsc,\cU_{s},\cU_{s,n})
\]
is relatively compact. By Lemma~\ref{lem:kappacompact}, it remains to see that
$X=\Un_{s\in\bbN^*}K_s$.

Assume that some element $x\in X$ is not in $\Un_{s\in\bbN^*}K_s$.
\be
\item Since $x\notin K_{()}$, there is $m_1$ with $x\notin \Un\sigma(\cU_{m_1})$.
\item Since $x\notin K_{m_1}$, there is $m_2$ with $x\notin \Un\sigma(\cU_{m_1,m_2})$.
\item Since $x\notin K_{m_1,m_2}$, there is $m_3$ with $x\notin \Un\sigma(\cU_{m_1,m_2,m_3})$.
\item Etc.
\ee
Then the play
\[
\cU_{m_1},\sigma(\cU_{m_1}),
\cU_{m_1,m_2}, \sigma(\cU_{m_1,m_2}),
\cU_{m_1,m_2,m_3}, \sigma(\cU_{m_1,m_2,m_3}),
\dotsc
\]
is lost by Bob; a contradiction.
\epf

Let $\alpha$ be an ordinal number.
The transfinite Menger game $\gfin^\alpha(\Op,\Op)$ is defined as the ordinary Menger game,
with the only difference that now there is an inning per each ordinal number $\beta<\alpha$.
Clearly, if $\alpha_1<\alpha_2$ and Bob has a winning strategy in the $\alpha_1$-Menger game, then
Bob has a winning strategy in the $\alpha_2$-Menger game: He can use a winning startegy in the first
$\alpha_1$ innings, and then play arbitrarily.
Thus, the following theorem is stronger than Theorem~\ref{thm:Mgame1}, and has no assumption on the topological space $X$.

The \emph{weight} of a topological space is the minimal cardinality of a base for its topology.

\bthm
\label{thm:Mgame2}
Let $X$ be a topological space of weight $\kappa$.
Bob has a winning strategy in the game $\gfin^\kappa(\Op,\Op)$ if and only if the space $X$ is a union of at most $\kappa$ compact sets.
\ethm
\bpf
$(\Impl)$ The proof is identical to that of Theorem~\ref{thm:Mgame1}, only that here we begin with a base of
cardinality $\kappa$.

$(\Leftarrow)$ In the $\alpha$-th inning, Bob covers the $\alpha$-th compact set.
\epf

\bcor
\label{cor:Mgame2}
Let $X$ be a topological space of weight $\kappa$.
If Bob has a winning strategy in the Menger game, then the space $X$ is a union of at most $\kappa$ compact sets.\qed
\ecor

The converse of Corollary~\ref{cor:Mgame2} is false: The discrete space of cardinality $\kappa$ has
weight $\kappa$, and it is a union of $\kappa$ compact sets (singletons).
This space is not Lindel\"of, and thus not Menger, so Bob has no winning strategy in the Menger game played on this space.

\medskip

For a space $X$, let $\cl{\Op}(X)$ be the families $\cU$ of open sets with $\Un_{U\in\cU}\cl{U}=X$.
A space $X$ is \emph{almost Menger} if it satisfies $\sfin(\Op,\cl{\Op})$.
For regular spaces, almost Menger is equivalent to Menger~\cite{BPS12}.
And similarly for the other notions considered below.
Thus, the remainder of this section is mainly relevant for nonregular spaces.
The \emph{almost Menger game} is the game associated to the property $\sfin(\Op,\cl{\Op})$.
The corresponding notion of \emph{almost Lindel\"of} is classic, and so is the notion of \emph{almost compact} space:
A space $K$ is almost compact if every open cover of $K$ has a finite subset $\cF$ with
dense union. This notion appears in the literature under various names. For Hausdorff spaces,
it is known to be equivalent to \emph{Hausdorff closed}, that is, being closed in all Hausdorff superspaces.

A set $K$ in a space $X$ is  \emph{relatively} almost compact if every open cover of the space $X$ has a finite subset $\cF$ with
$K\sub \Un _{U\in\cF}\cl{U}$. The standard proof that every almost compact space is $H$-closed shows that
every relatively almost compact set in a Hausdorff space is closed in that space.
However, no separation hypothesis is needed in the following theorem.
Since the proof of the following theorem is provided by the first part of Scheepers's argument,
we attribute the theorem to Scheepers.

\bthm[Scheepers]
\label{thm:aMgame}
Let $X$ be a second countable topological space. The following assertions are equivalent:
\be
\item Bob has a winning strategy in the game $\gfin(\Op,\cl{\Op})$
played on $X$.
\item The space $X$ is a countable union of relatively almost compact sets.
\ee
\ethm
\bpf
$(2)\Impl (1)$: This is easy.

$(1)\Impl (2)$:
This is a part of the argument of Scheepers~\cite[Theorem~1]{Sch95}.
We provide it, for completion and verification.

Repeat the inductive construction of Theorem~\ref{thm:Mgame1} verbatim.
Having completed it, define for each sequence $s\in\bbN^*$:
\[
K_s := \CapN{n}\Un_{U\in\sigma(\cU_{s_1},\dotsc,\cU_{s},\cU_{s,n})}\cl{U}.
\]
The set $K_s$ is relatively almost compact.
It remains to see that $X=\Un_{s\in\bbN^*}K_s$.

Assume that some element $x\in X$ is not in $\Un_{s\in\bbN^*}K_s$.
\be
\item Since $x\notin K_{()}$, there is $m_1$ with $x\notin \Un_{U\in\sigma(\cU_{m_1})}\cl{U}$.
\item Since $x\notin K_{m_1}$, there is $m_2$ with $x\notin \Un_{U\in\sigma(\cU_{m_1,m_2})}\cl{U}$.
\item Since $x\notin K_{m_1,m_2}$, there is $m_3$ with $x\notin \Un_{U\in\sigma(\cU_{m_1,m_2,m_3})}\cl{U}$.
\item Etc.
\ee
Then the play
\[
\cU_{m_1},\sigma(\cU_{m_1}),
\cU_{m_1,m_2}, \sigma(\cU_{m_1,m_2}),
\cU_{m_1,m_2,m_3}, \sigma(\cU_{m_1,m_2,m_3}),
\dotsc
\]
is lost by Bob; a contradiction.
\epf

The assertions analogous to the more general Theorem~\ref{thm:Mgame2} and Corollary~\ref{cor:Mgame2} also hold.
Theorem~\ref{thm:aMgame} answers a question of Babinkostova, Pansera and Scheepers~\cite[Question~26(2)]{BPS12}
in the case of second countable spaces.

\section{The additivity of Rothberger's property}
\label{sec:add}

Let $\add(\cN)$ be the minimal cardinality of a family $F \sub\NN$
such that there is no function $S\colon \bbN\to\Fin$ with $\card{S(n)} \leq n$ for all $n$,
such that  for each function $f \in F$ we have $f(n) \in S(n)$ for all but finitely many $n$.

The notation $\add(\cN)$ is explained by a result of Bartoszy\'nski and Judah~\cite[Theorem~2.11]{covM2}:
The cardinal number $\add(\cN)$ is the minimal cardinality of a family of Lebesgue null sets of real numbers whose union is not Lebesgue null.
In general, the \emph{additivity} of a property is the minimum cardinality of a family of sets satisfying the property,
whose union does not.
The following theorem is attributed to Carlson by Bartoszy\'nski and Judah~\cite[Theorem~2.9]{covM2}.

\bthm[Carlson]
\label{thm:Carlson}
Let $\kappa<\add(\cN)$. If a Lindel\"of space is a union of at most $\kappa$ spaces satisfying $\sone(\Op,\Op)$, then the space $X$ satisfies $\sone(\Op,\Op)$.
That is, for Lindel\"of spaces, $\add(\cN)\le \add(\sone(\Op,\Op))$.
\ethm

This theorem is an easy consequence of a simple, basic fact concerning selection principles.
We need the following lemmata.

%Essentially 2.3.9 in the book

Let $\rmA$ and $\rmB$ be types of open covers.
A topological space $X$ satisfies $\mathsf{S}_{n}(\rmA,\rmB)$ if for all $\cU_1,\cU_2,\dots\in\rmA(X)$,
there are finite sets $\cF_1\sub\cU_1,\cF_2\sub\cU_2,\dotsc$ such that
$\card{\cF_n}\le n$ for all $n$, and $\UnN{n}\cF_n\in\rmB(X)$.

Garcia-Ferreira and Tamariz-Mascarua~\cite[Lemma 3.12]{GFTM95} established the following observation
in the case $\rmA=\Op$.

\blem[{\cite[Theorem~A.1]{MHP}}]
\label{lem:sn}
Let $\rmA$ be a type of countable covers such that
every pair of covers of type $\rmA$ has a joint refinement of type $\rmA$.
Then $\mathsf{S}_{n}(\rmA,\Op)=\sone(\rmA,\Op)$.
\elem

\blem[Folklore]
If a space $X$ satisfies $\sone(\rmA,\Op)$, then
for each sequence $\eseq{\cU}$ in $\Op(X)$ there are elements
$U_{m_1}\in\cU_{m_1},U_{m_2}\in\cU_{m_2},\dotsc$ such that
for each point $x\in X$, we have $x\in U_{m_n}$ for infinitely many $n$.
\elem
\bpf
As usual, we split the sequence of open covers to infinitely many disjoint sequences, and apply
the property $\sone(\rmA,\Op)$ to each subsequence separately.
\epf

\bthm
\label{thm:carlsongen}
Let $\rmA$ be a type of countable covers such that
%every cover of type $\rmA$ has a countable refinement of type $\rmA$
%and
every pair of covers of type $\rmA$ has a joint refinement of type $\rmA$.
Then $\add(\cN)\le\add(\sone(\rmA,\Op))$.
\ethm
\bpf
Let $\kappa<\add(\cN)$ and $X = \Un_{\alpha < \kappa}X_\alpha$, where each space
$X_\alpha$ satisfies $\sone(\rmA,\Op)$.
By Lemma~\ref{lem:sn}, it suffices to show that the space $X$ satisfies
$\mathsf{S}_{n}(\rmA,\Op)$.

Let $\cU_n=\{U^n_m : m\in\bbN\}\in\rmA(X)$, for $n\in\bbN$.
For each ordinal number $\alpha<\kappa$, as the space $X_\alpha$ satisfies $\sone(\rmA,\Op)$,
there is a function $f_\alpha\in\NN$ such that
for each point $x\in X_\alpha$ we have
\[
x\in U^n_{f_\alpha(n)}
\]
for infinitely many $n$.

There is a function
$S\colon \bbN\to \Fin$ with $\card{S(n)} \leq n$ for all $n$,
such that for each $\alpha<\kappa$ we have
\[
\ f_\alpha(n) \in S(n)
\]
for all but finitely many $n$.
Then $\UnN{n}\sset{U^n_{m}}{m\in S(n)}\in\Op(X)$.
\epf

Judging by an extensive survey on the topic~\cite{AddQuad},
the result in the second item below seems to be new.

\bcor
\mbox{}
\be
\item For Lindel\"of spaces, $\add(\cN)\le\add(\sone(\Op,\Op))$.
\item $\add(\cN)\le\add(\sone(\Ga,\Op))$.
\ee
\ecor
\bpf
(1) Since the spaces are Lindel\"of, we may restrict attention to countable covers,
and the assumptions of Theorem~\ref{thm:carlsongen} hold.

(2) A countably infinite subset of a point-cofinite cover is also a point-cofinite cover. Thus, we may
restrict attention to countable point-cofinite covers.
It is well-known that every pair of point-cofinite covers has a joint refinement that is a point-cofinite cover.
Indeed, let $\cU$ and $\cV$ be countable point-cofinite covers.
Enumerate them $\cU=\sseq{U_n}$ and $\cV=\sseq{V_n}$.
Then $\sseq{U_n\cap V_n}\in\Ga(X)$. Theorem~\ref{thm:carlsongen} applies.
\epf

We can extract additional information from this proof method.
A topological space $X$ satisfies $\mathsf{U}_n(\Ga,\Ga)$~\cite{MHP} if
for all $\cU_1,\cU_2,\dots\in\Ga(X)$,
there are finite sets $\cF_1\sub\cU_1,\cF_2\sub\cU_2,\dotsc$ such that
$\card{\cF_n}\le n$ for all $n$, and $\sseq{\Un\cF_n}\in\Ga(X)$.
This property is strictly inbetween $\sone(\Ga,\Ga)$ and $\ufin(\Op,\Ga)$~\cite[Theorems~3.3 and 3.8]{MHP}.

\bthm
$\add(\cN)\le\add(\mathsf{U}_n(\Ga,\Ga))$.
\ethm
\bpf
Let $\kappa<\add(\cN)$ and $X = \Un_{\alpha < \kappa}X_\alpha$, where each space
$X_\alpha$ satisfies $\mathsf{U}_n(\Ga,\Ga)$.
It suffices to show that the space $X$ satisfies
$\mathsf{U}_{n^2}(\Ga,\Ga)$, where the cardinality of the $n$-th selected finite set is at most $n^2$~\cite[Lemma~3.2]{MHP}.

Let $\cU_n=\{U^n_m : m\in\bbN\}\in\Ga(X)$, for $n\in\bbN$.
For each $\alpha<\kappa$, as the space $X_\alpha$ satisfies $\mathsf{U}_n(\Ga,\Ga)$,
there is a function $S_\alpha\colon\bbN\to \prod_n[\bbN]^{\le n}$ such that
for each point $x\in X_\alpha$ we have
\[
x\in \Un_{m\in S_\alpha(n)}U^n_m
\]
for infinitely many $n$.

There is a function
$S\colon \bbN\to \prod_n[\bbN]^{\le n}$ with $\card{S(n)} \leq n$ for all $n$,
such that for each $\alpha<\kappa$ we have
\[
\ S_\alpha(n) \in S(n)
\]
for all but finitely many $n$.
For each natural number $n$, let $F_n:=\Un S(n)$.
Then $\card{F_n}\le n^2$ for all $n$, and
$\sseq{\Un_{m\in F_n}U^n_m}\in\Ga(X)$.
\epf

\subsection*{Acknowledgments}
This is the first paper that I write after a long, callenging period.
I thank all those who helped and encouraged me throughout, and supported my
return to normal track afterwards.
Above all, I thank my wife, Adina, for her faith, support, and patience.

\ed